\mag 1200
\input amstex
\input amsppt.sty

\overfullrule 0pt

\let\al\alpha
\let\bt\beta
\let\gm\gamma 
\let\dl\delta 
 \let\eps\varepsilon \let\epsilon\eps

\let\la\lambda \let\La\Lambda
 
 \let\phi\varphi

\def\C{\Bbb C}

\def\Z{\Bbb Z}

\def\F{\Cal F}

  \def\Fwh{\,\wh{\!\F}}

\def\A{\roman A}

\def\d{\roman d}

\def\E{\roman E}

\let\Fun\Fwh


\def\A{\Cal A} 

\def\Fun/{ {\text{Fun}}} 
\def\End/{ {\text{End}}} 
\def\Hom/{ {\text{Hom}}} 
\def\Rm/{\^{$R$-}matrix}

\def\M/{\Cal M } 
\def\A/{$A_{\tau,\eta}(sl_2)$} 
\def\Am/{$A_{\tau,\eta,\mu}(sl_2)$} 

\def\E/{$E_{\tau,\eta}(sl_2)$} 

\def\Ex/{ \widetilde {  E   }}
\def\EE/{ $\widetilde {  E   }_{\tau,\eta}(sl_2)$}

\def\eqg/{elliptic quantum group}
\def\Em/{$E_{\tau,\eta,\mu}(sl_2)$} 
\def\E/{$E_{\tau,\eta}(sl_2)$} 

\def\a{ \widetilde a}
\def\b{ \widetilde b}
\def\c{ \widetilde c}
\def\d{ \widetilde d}
\def\h{ \widetilde h}
\def\t{ \widetilde t}

\def\Det{ {\text{Det}}}

\def\Dett{ \widetilde {\text{Det}}}

\def\Ce/{\C/{1\over \eta}\Z }
\def\V{\Cal V }

\topmatter

\title
On representations of the elliptic quantum group
$E_{\tau,\eta}(sl_2)$
\endtitle

\author
Giovanni Felder and Alexander Varchenko
\endauthor

\thanks
The authors were supported in part by NSF grants
DMS-9400841 and DMS-9501290
\endthanks

\date
December, 1995
\enddate

\abstract
We describe representation theory of the \eqg/ \E/. It turns out 
that the representation theory  is parallel to 
the representation
theory of the Yangian $Y(sl_2)$ and the  quantum loop group
$ U_q(\widetilde {  sl   }_2)$.

We introduce basic notions of 
representation theory of the elliptic quantum group
\E/ and construct three families of modules: evaluation modules, cyclic
modules, one-dimensional modules. 
We show that under certain conditions any irreducible highest weight
module of finite type is isomorphic to a tensor product of evaluation modules
and a one-dimensional module. 
We describe fusion of finite dimensional evaluation modules. 
In particular, we show that under certain conditions
the tensor product of two evaluation modules becomes
reducible and contains an
 evaluation module, in this case the imbedding 
of the evaluation module into the tensor product
 is given in terms of elliptic binomial coefficients.
We describe the determinant element of the \eqg/.
Representation theory becomes special if $N\eta = m + \ell \tau$ where
$N, m, \ell$ are integers. We indicate some new features in this case.
\endabstract

 \leftheadtext{ G.Felder and A.Varchenko }
  \rightheadtext{ On representations of the elliptic quantum group}
\endtopmatter

\document

\head 1. Introduction \endhead

The elliptic quantum group is an algebraic structure underlying
the elliptic solutions of the Star-Triangle relation in statistical
mechanics and connected with the Knizhnik-Zamolodchikov-Bernard
equation on torus. In this paper we consider the
\eqg/ \E/ defined in [Fe1-2] and discuss its representation theory.
It turns out that representation theory 
of the \eqg/ \E/ is parallel to representation
theory of the Yangian $Y(sl_2)$ and the quantum loop group
$ U_q(\widetilde {  sl   }_2)$, cf. [CP].

We  introduce 
basic notions of representation theory
of the elliptic quantum group (notions of the operator algebra,
 a highest weight module,
an irreducible
module, a module of finite type, a
singular vector). Essentially  all the
 notions are formulated in terms of the
associated operator algebra.

We construct three families of \E/-modules: evaluation Verma
modules $V_\La(z)$,
 cyclic
modules, and one-dimensional modules.

An evaluation Verma module
is an infinite dimensional module
determined by an evaluation point $z\in \C$
 and  a weight $\La \in \C$. If $\La = n + (m + \ell \tau)/2\eta$,
where $n, m, \ell$ are integers, $n \geq 0$, 
then the Verma module has a submodule and the quotient module $L_\La(z)$
has
finite dimension $n+1$. For generic $\eta$ the evaluation Verma module
$V_\La(z)$
is irreducible unless $\La = n + (m + \ell \tau)/2\eta$. If
$\La$ has this form,  then $L_\La(z)$ is irreducible.

A tensor product of irreducible evaluation modules is an irreducible
highest weight module of finite type, if $\eta$ is generic and 
the evaluation
points of factors
are generic. We show that under certain conditions  every irreducible
highest weight module of finite type is isomorphic to a tensor product
of irreducible evaluation modules and a one-dimensional module.

 We give necessary and sufficient
conditions for the  tensor product of two finite dimensional
modules $L_{\La_1}(z_1) \otimes L_{\La_2}(z_2) $
to be reducible. In particular, 
$L_{\La_1}(z_1) \otimes L_{\La_2}(z_2) $  becomes reducible, if
$z_2-z_1=(\La_1+\La_2)\eta$. In this case the tensor product contains
a submodules isomorphic to $L_{\La_1+\La_2}(z_2-\La_1\eta)$. The imbedding
of this module into the tensor product is given in terms of  elliptic
binomial coefficients.

We indicate the determinant element of the \eqg/. It is a group-like
central element, see precise statements in Section 10. 

The elliptic quantum group depends on complex parameters $\tau$ and
$\eta$. We show that the \eqg/s with parameters $(\tau, \eta)$,
$(\tau + 1, \eta)$, and $(-1/\tau, -\eta/\tau)$ are essentially isomorphic.

Representation theory of the elliptic quantum group is special,
if $N\eta = m + \ell \tau$ where $N, m, \ell$ are integers. We indicate some
new features in this case.

Now we
briefly sketch the definition of the \eqg/ \E/  and will give a detailed
account in Sections 2 and 3.

To define the \eqg/ \E/ we fix two nonzero
complex numbers $\tau$ and $\eta$, Im $ \tau > 0$,  and start from a
$4 \times 4$-matrix $R(\la,w) \in \End/(\C^2\otimes \C^2)$,
$$
\gather
R(\la,w,\eta,\tau) = E_{1,1}\otimes E_{1,1}+E_{2,2}\otimes E_{2,2}+
\al(\la,w,\eta,\tau) E_{1,1}\otimes E_{2,2} + \\
\beta(\la,w,\eta,\tau)
E_{1,2}\otimes E_{2,1}
+ \gm(\la,w,\eta,\tau) E_{2,1}\otimes E_{1,2} +
\dl(\la,w,\eta,\tau) E_{2,2}\otimes E_{1,1}   \\
\endgather
$$
where the functions $\al$, $\beta$, $\gm$,
$\dl$ are defined in Section 2, $E_{i,j}$
is the $2\times 2$-matrix with the only nonzero element $1$ at the intersection
of the $i$-th row and $j$-th column.

Let $h$ be the diagonal $2\times 2$ matrix ${\text{Diag}}(1,-1)$.
The $R$-matrix
 $R(w,\la)$ satisfies the modified
Yang-Baxter equation  in $\End/(\C^2\otimes \C^2\otimes \C^2)$,
$$
\gather
R^{(12)}(\la-2\eta h^{(3)},w_{12})\, R^{(13)}(\la,w_{13})\,
R^{(23)}(\la-2\eta h^{(1)},w_{23})\, =           \\
R^{(23)}(\la,w_{23})\,R^{(13)}(\la-2\eta h^{(2)},w_{13})\,
R^{(12)}(\la,w_{12})\, .\\
\endgather
$$
Here $w_{ij}=w_i-w_j$;
$R^{(12)}(\la-2\eta h^{(3)},w_{12})$ means that if
$a\otimes b\otimes c \in\C^2\otimes \C^2\otimes \C^2$ and $h c = \mu c$,
 $\mu \in \C$,
then $R^{(12)}(\la-2\eta h^{(3)},w_{12})\, a\otimes b\otimes c =
 R^{(12)}(\la-2\eta \mu ,w_{12})( a\otimes b) \otimes c $, and the other
symbols have a similar meaning.

The \eqg/ \E/ is an algebra generated by meromorphic functions of
a variable $h$ and the matrix elements of a matrix $L(\la,w) \in
\End/(\C^2)$ with non commutative entries, subject to the relations
$$
\gather
R^{(12)}(\la-2\eta h, w_{12})\, L^{(1)}(\la, w_{1})\,
L^{(2)}(\la-2\eta h^{(1)}, w_{2})\, =   \tag 1        \\
L^{(2)}(\la, w_{2})\,L^{(1)}(\la-2\eta h^{(2)}, w_{1})\,
R^{(12)}(\la, w_{12})\, .            \\
\endgather
$$
Here $h$ is considered as a generator of a one dimensional commutative Lie
algebra $\frak h$.

An \E/-module is a 
diagonalizable $\frak h$-module $V$ together with a meromorphic
function $L(\la,w)$ on $\frak h \times \C$ with values in
$\End/(\C^2 \otimes V)$ such that identity (1) holds in
$\End/(\C^2 \otimes \C^2 \otimes V)$ and so that $L$ is of weight zero:
$$
[h^{(1)} + h^{(2)}\, , \, L(\la,w)]\,=\,0.
$$

If $V$ and $U$ are \E/-modules, then $V\otimes U$ is an \E/-module
with an $\frak h$-module structure $h\,(a\otimes b)\,=\,
h a\otimes b+a\otimes h b$ and an $L$-operator
$$
L^{(12)}(\la-2\eta h^{(3)},w) \,  L^{(13)}(\la,w).
$$
If $V,U,W$ are \E/-modules, then the modules
$(V\otimes U)\otimes W$ and $V\otimes (U\otimes W)$
are isomorphic with the obvious isomorphism, see [Fe1-2].

Similarly one could define the elliptic quantum group
associated with a
simple Lie algebra of type $A, B, C, D$, see [Fe2].

\head 2. The elliptic quantum group \E/ \endhead

Let 
$$
\theta(z,\tau) = - \sum^\infty_{j=-\infty} e^{\pi i(j+{1\over 2})^2\tau
+2\pi i (j+ {1\over 2})(z+{1\over 2})}
$$
be the Jacobi theta function, and let functions $\al, 
\bt, 
\gm,
\dl$ be given by
$$
\al(w,\la,\eta,\tau)={\theta(w)\theta(\la+2\eta) \over 
\theta(w-2\eta)\theta(\la) }              
\qquad 
\bt(w,\la,\eta,\tau)={\theta(-w-\la)\theta(2\eta) \over 
\theta(w-2\eta)\theta(\la) }       
$$
$$
\gm(w,\la,\eta,\tau)={\theta(w-\la)\theta(2\eta) \over 
\theta(w-2\eta)\theta(\la) }              
\qquad
\dl(w,\la,\eta,\tau)={\theta(w)\theta(\la-2\eta) \over 
\theta(w-2\eta)\theta(\la) }                     
$$

Let $\eta$ and $\tau$ be nonzero complex numbers, Im $\tau >  0$.

The \eqg/ \E/ is the algebra over $\C$
with generators of two types. The generators of the first type
are labelled by meromorphic functions $f(h)$
of one complex variable with period $1/\eta$,
$f(h+1/\eta)=f(h)$.
The generators of the second type are
$a(\la,w),\, b(\la,w),
c(\la,w),\, d(\la,w)$. They  are labelled by elements $\la \in
\C/\Z$ and complex numbers $w \in \C$.
 The generators of the
algebra satisfy the following two groups of relations .
       
Relations involving $h$:
$$
\gather
f(h)\,g(h) \, = \,g(h)\,f(h) , \tag 2  \\
f(h)\,a(\la,w)\,=\,  a(\la,w) \, f(h) , \qquad
f(h)\,d(\la,w)\,=\,  d(\la,w) \, f(h)  , \\
f(h)\,b(\la,w)\,=\,  b(\la,w) \, f(h-2) , \qquad
f(h)\,c(\la,w)\,=\,  c(\la,w) \, f(h+2)  , \\
\endgather
$$
where $f(h), g(h)$ are  generators of the first type. 

The remaining relations have the form:
$$
\gather                %
a_1(\la)a_2(\la\!-\!2\eta) %
=a_2(\la)a_1(\la\!-\!2\eta), \tag 3 \\     %
a_1(\la)b_2(\la\!-\!2\eta) %
=b_2(\la)a_1(\la\!+\!2\eta)\alpha(\la)+
a_2(\la)b_1(\la\!-\!2\eta)\gamma(\la), \\     %
b_1(\la)a_2(\la\!+\!2\eta) %
=b_2(\la)a_1(\la\!+\!2\eta)\beta(\la)
+a_2(\la)b_1(\la\!-\!2\eta)\delta(\la), \\     %
b_1(\la)b_2(\la\!+\!2\eta) %
=b_2(\la)b_1(\la\!+\!2\eta), \\     %
\beta(\la\!-\!2\eta h)c_1(\la)a_2(\la\!-\!2\eta)
+\alpha(\la\!-\!2\eta h)a_1(\la)
c_2(\la\!-\!2\eta) %
=c_2(\la)a_1(\la\!-\!2\eta), \\     %
\beta(\la\!-\!2\eta h)c_1(\la)b_2(\la\!-\!2\eta)+
\alpha(\la\!-\!2\eta h)a_1(\la)
d_2(\la\!-\!2\eta) %
=d_2(\la)a_1(\la\!+\!2\eta)\alpha(\la)+
c_2(\la)b_1(\la\!-\!2\eta)\gamma(\la), \\     %
\beta(\la\!-\!2\eta h)d_1(\la)a_2(\la\!+
\!2\eta)+\alpha(\la\!-\!2\eta h)b_1(\la)
c_2(\la\!+\!2\eta) %
=d_2(\la)a_1(\la\!+\!2\eta)\beta(\la)+
c_2(\la)b_1(\la\!-\!2\eta)\delta(\la), \\     %
\beta(\la\!-\!2\eta h)d_1(\la)b_2(\la\!+
\!2\eta)+\alpha(\la\!-\!2\eta h)b_1(\la)
d_2(\la\!+\!2\eta) %
=d_2(\la)b_1(\la\!+\!2\eta), \\     %
\delta(\la\!-\!2\eta h)c_1(\la)a_2(\la\!-\!2\eta)+
\gamma(\la\!-\!2\eta h)a_1(\la)
c_2(\la\!-\!2\eta) %
=a_2(\la)c_1(\la\!-\!2\eta), \\     %
\delta(\la\!-\!2\eta h)c_1(\la)b_2(\la\!-\!2\eta)+
\gamma(\la\!-\!2\eta h)a_1(\la)d_2(\la\
\!-\!2\eta) %
=b_2(\la)c_1(\la\!+\!2\eta)\alpha(\la)+
a_2(\la)d_1(\la\!-\!2\eta)\gamma(\la), \\     %
\delta(\la\!-\!2\eta h)d_1(\la)a_2(\la\!+
\!2\eta)+\gamma(\la\!-\!2\eta h)b_1(\la)
c_2(\la\!+\!2\eta) %
=b_2(\la)c_1(\la\!+\!2\eta)\beta(\la)+
a_2(\la)d_1(\la\!-\!2\eta)\delta(\la), \\     %
\delta(\la\!-\!2\eta h)d_1(\la)b_2(\la\!+
\!2\eta)+\gamma(\la\!-\!2\eta h)b_1(\la)
d_2(\la\!+\!2\eta) %
=b_2(\la)d_1(\la\!+\!2\eta), \\     %
c_1(\la)c_2(\la\!-\!2\eta) %
=c_2(\la)c_1(\la\!-\!2\eta), \\     %
c_1(\la)d_2(\la\!-\!2\eta) %
=d_2(\la)c_1(\la\!+\!2\eta)\alpha(\la)+
c_2(\la)d_1(\la\!-\!2\eta)\gamma(\la), \\     %
d_1(\la)c_2(\la\!+\!2\eta) %
=d_2(\la)c_1(\la\!+\!2\eta)\beta(\la)+
c_2(\la)d_1(\la\!-\!2\eta)\delta(\la), \\     %
d_1(\la)d_2(\la\!+\!2\eta) %
=d_2(\la)d_1(\la\!+\!2\eta)     %
\endgather           %
$$
where $\al, \bt, \gm, \dl$ are functions of $w=w_1-w_2, h,
 \la, \eta, \tau$,  the symbol $a_j(\la)$ denotes
the generator $a(\la,w_j)$, and similarly for $b_j, c_j, d_j$.

Note, for example, that $\al(w_1-w_2,\la- 2\eta h,\eta,\tau)$
as a function of $h$ is meromorphic and $1/\eta$-periodic, therefore, it
is an element of the \eqg/.

The functions $\al, \beta, \gm, \dl$ have a pole at $w_1-w_2=2\eta$. In this 
special case the relations take the form
$$
\gather                %
a(\la,w+2\eta)a(\la\!-\!2\eta,w) %
=a(\la,w)a(\la\!-\!2\eta,w+2\eta),    \tag 4 \\     %
b(\la,w+2\eta)b(\la\!+\!2\eta,w) %
=b(\la,w)b(\la\!+\!2\eta,w+2\eta), \\     %
c(\la,w+2\eta)c(\la\!-\!2\eta,w) %
=c(\la,w)c(\la\!-\!2\eta,w+2\eta), \\     %
d(\la,w+2\eta)d(\la\!+\!2\eta,w) %
=d(\la,w)d(\la\!+\!2\eta,w+2\eta),  \\
c(\la,w+2\eta)a(\la\!-\!2\eta,w)=a(\la,w+2\eta)
c(\la\!-\!2\eta,w), \\     %
d(\la,w+2\eta)b(\la\!+\!2\eta,w)=b(\la,w+2\eta)
d(\la\!+\!2\eta,w) , \\     %
\theta(\la+2\eta)b(\la,w)a(\la\!+\!2\eta,w+2\eta)
=\theta(\la-2\eta)a(\la,w)b(\la\!-\!2\eta,w+2\eta), \\     %
\theta(\la+2\eta)d(\la,w)c(\la\!+\!2\eta,w+2\eta)=
\theta(\la-2\eta)c(\la,w)d(\la\!-\!2\eta,w+2\eta), \\     %
\endgather           
$$
$$
\gather
 d(\la,w+2\eta)a(\la+2\eta,w)
- b(\la,w+2\eta)c(\la+2\eta,w) \,=   \\
=\, a(\la,w+2\eta)d(\la-2\eta,w)
- c(\la,w+2\eta)b(\la-2\eta,w)\,= \\
\endgather
$$
$$
\gather
{\theta(\la-2\eta h-2\eta)\over \theta(\la-2\eta h)} 
( d(\la,w+2\eta)a(\la+2\eta,w)
- b(\la,w+2\eta)c(\la+2\eta,w)) \,=   \\
=\,{\theta(\la-2\eta)\over \theta(\la)}
a(\la,w)
d(\la\!-\!2\eta,w+2\eta) - {\theta(\la+2\eta)\over \theta(\la)}
b(\la,w)c(\la\!+\!2\eta,w+2\eta), \\
\endgather
$$
$$
\gather
{\theta(\la-2\eta h+2\eta)\over \theta(\la-2\eta h)} 
( a(\la,w+2\eta)d(\la-2\eta,w)\,
-\, c(\la,w+2\eta)b(\la-2\eta,w))\,= \\
=\,{\theta(\la+2\eta)\over \theta(\la)}
d(\la,w)a(\la\!+\!2\eta,w+2\eta)\, -\,
{\theta(\la-2\eta)\over \theta(\la)}
c(\la,w)b(\la\!-\!2\eta,w+2\eta)\, . \\
\endgather
$$
These relations are obvious regularizations of relations (3).
Similar relations hold for other poles of the functions
$\al,\,\beta,\,\gm,\,\delta .$

An  {\it \E/-module structure} on a complex vector space $V$ is 
a direct sum decomposition 
$$
V\,=\, \oplus_{\mu \in \C/{1\over \eta}\Z} V[\mu]
$$
and endomorphisms 
$ a(\la,w),$ $ b(\la,w),$
$c(\la,w),$  $ d(\la,w) \in \End/(V)$ which meromorphically depend on
$\la \in \C/\Z$, 
$ w\in \C$. The direct sum decomposition allows us to define
endomorphisms $f(h) \in \End/(V)$ by the rule $f(h)v=f(\mu)v$, if
$v\in V[\mu]$. We assume that the endomorphisms $f(h)$,
$ a(\la,w),$ $ b(\la,w),$
$c(\la,w),$  $ d(\la,w)$ satisfy the relations of the \eqg/.

 We allow the fact
that for a given module structure the action of 
 not
all elements of the \eqg/ is well defined.

A module is of {\it finite type} if each 
space $V[\mu]$ has finite dimension.

If $V, \, W$ are \E/-modules, then the tensor product $V\otimes W$ has
an \E/-module structure, where the action of 
$f(h),\, a(\la,w),\, b(\la,w),
c(\la,w),\, d(\la,w)$ is given by
$$
\gather
f(h^{(1)}+h^{(2)}) ,\, \\
a(\la-2\eta h^{(2)},w)\otimes a(\la,w) \, + \, 
b(\la-2\eta h^{(2)},w)\otimes c(\la,w) ,\\
a(\la-2\eta h^{(2)},w)\otimes b(\la,w) \, + \, 
b(\la-2\eta h^{(2)},w)\otimes d(\la,w) ,\\
c(\la-2\eta h^{(2)},w)\otimes a(\la,w) \, + \, 
d(\la-2\eta h^{(2)},w)\otimes c(\la,w) ,\\
c(\la-2\eta h^{(2)},w)\otimes b(\la,w) \, + \, 
d(\la-2\eta h^{(2)},w)\otimes d(\la,w) ,\\
\endgather
$$
here $a(\la-2\eta h^{(2)},w)\otimes b(\la,w) $ 
means that if $s\otimes t \in V\otimes W$
and $b(\la,w)t \in W[\mu]$, then
$a(\la-2\eta h^{(2)},w) \otimes b(\la,w) \, s\otimes t \, = \,
a(\la-2\eta \mu, w)\,s\,\otimes b(\la,w) \,t $; if 
$s\otimes t \in V[\nu]\otimes W[\mu]$ then $f(h^{(1)} + h^{(2)})$ 
$s\otimes t$ $= $ $f(\nu+\mu)s\otimes t$;
and the other terms have similar 
meaning.

The tensor product is associative.

\head 3. The operator algebra  \endhead

In this section we introduce the 
operator algebra of the elliptic quantum group,
cf. [ABB], [BBB].

For a complex vector space $V$ let
 $\Fun/(V)$ be the space of meromorphic functions of $\la$ with values
in $V$ and $1$-periodic, $F(\la+1)=F(\la)$.
  The space $\Fun/(V)$ has a natural structure of a vector space
over the field  $\Fun/(\C)$.

For an \E/-module $V$ define its {\it operator algebra} $A(V)$
as the algebra of the following  
operators
 $ f(\h), \a(w), \b(w), \c(w), \d(w)$ acting on  $\Fun/(V)$.
 Introduce endomorphisms $ f(\h), \, \a(w), \, \b(w), \, \c(w), \,
\d(w) \in \End/(\Fun/(V))$ by the rule
$$
\gather
(f(\h) F) (\la)\,=\, f(h) \, F(\la), \\
(\a(w) F) (\la)\,=\, a(\la,w) \, F(\la-2\eta), \qquad 
(\c(w) F) (\la)\,=\, c(\la,w) \, F(\la-2\eta), \\
(\b(w) F) (\la)\,=\, b(\la,w) \, F(\la+2\eta), \qquad 
(\d(w) F) (\la)\,=\, d(\la,w) \, F(\la+2\eta), \\
\endgather
$$
where $F \in \Fun/(V)$. 

 The relations for the generators in the \eqg/ induce some
universal relations for operators
$ f(\h), \, \a(w), \, \b(w), \, \c(w), \, \d(w)$
and motivate the following definition of the {\it operator algebra}
\A/ of the \eqg/. The operator algebra \A/ is an algebra over $\C$
 generated by $ f(\la,\h), \, \a(w), $ $ \b(w),$ $ \c(w),$ $\d(w)$,
where $f(\la,\h)$ runs through the space of meromorphic functions
of two variables $\la$ and $\h$ which are $1$-periodic 
in $\la$ and
$1/\eta$-periodic 
in $\h$,
$$
f(\la+1,\h)\,=\,f(\la,\h),\qquad f(\la,\h+1/\eta)\,=\,f(\la,\h).
$$
 The generators of the algebra
satisfy the following two groups of relations.

Relations involving $f(\la,\h)$:
$$
\gather
f(\la,\h)\,g(\la,\h) \, = \,g(\la,\h)\,f(\la,\h) , \\
f(\la-2\eta,\h)\,\a(w)\,=\,  \a(w) \, f(\la,\h) , \qquad
f(\la +2\eta,\h)\,\d(w)\,=\,  \d(w) \, f(\la,\h)  , \\
f(\la+2\eta,\h+2)\,\b(w)\,=\,  \b(w) \, f(\la, \h) , \qquad
f(\la-2\eta,\h-2)\,\c(w)\,=\,  \c(w) \, f(\la,\h)  . \\
\endgather
$$

The remaining relations have the form:
$$
\gather                %
\a(w_1)\a(w_2) %
=\a(w_2)\a(w_1), \\     %
\a(w_1)\b(w_2) %
=\alpha(\la)\b(w_2)\a(w_1)+\gamma(\la)\a(w_2)\b(w_1), \\     %
\b(w_1)\a(w_2) %
=\beta(\la)\b(w_2)\a(w_1)+\delta(\la)\a(w_2)\b(w_1), \\     %
\b(w_1)\b(w_2) %
=\b(w_2)\b(w_1), \\     %
\beta(\la\!-\!2\eta \h)c(w_1)\a(w_2)+\alpha(\la\!-\!2\eta \h)\a(w_1)
\c(w_2) %
=\c(w_2)\a(w_1), \\     %
\beta(\la\!-\!2\eta \h)\c(w_1)\b_(w_2)+\alpha(\la\!-\!2\eta \h)\a(w_1)
\d(w_2) %
=\alpha(\la)\d(w_2)\a(w_1)+\gamma(\la)\c(w_2)\b(w_1), \\     %
\beta(\la\!-\!2\eta \h)\d(w_1)\a(w_2)+\alpha(\la\!-\!2\eta \h)\b(w_1)
\c(w_2) %
=\beta(\la)\d(w_2)\a(w_1)+\delta(\la)\c(w_2)\b(w_1), \\     %
\beta(\la\!-\!2\eta \h)\d(w_1)\b(w_2)+\alpha(\la\!-\!2\eta \h)\b(w_1)
\d(w_2) %
=\d(w_2)\b(w_1), \\     %
\delta(\la\!-\!2\eta \h)\c(w_1)\a(w_2)+\gamma(\la\!-\!2\eta \h)\a(w_1)
\c(w_2) %
=\a(w_2)\c(w_1), \\     %
\delta(\la\!-\!2\eta \h)
\c(w_1)\b(w_2)+\gamma(\la\!-\!2\eta \h)\a(w_1)\d(w_2)  %
=\alpha(\la)\b(w_2)\c(w_1)+\gamma(\la)\a(w_2)\d(w_1), \\     %
\delta(\la\!-\!2\eta \h)\d(w_1)\a(w_2)+\gamma(\la\!-\!2\eta \h)\b(w_1)
\c(w_2) %
=\beta(\la)\b(w_2)\c(w_1)+\delta(\la)\a(w_2)\d(w_1), \\     %
\delta(\la\!-\!2\eta \h)\d(w_1)\b(w_2)+\gamma(\la\!-\!2\eta \h)\b(w_1)
\d(w_2) %
=\b(w_2)\d(w_1), \\     %
\c(w_1)\c(w_2) %
=\c(w_2)\c(w_1), \\     %
\c(w_1)\d(w_2) %
=\alpha(\la)\d(w_2)\c(w_1)+\gamma(\la)\c(w_2)\d(w_1), \\     %
\d(w_1)\c(w_2) %
=\beta(\la)\d(w_2)\c(w_1)+\delta(\la)\c(w_2)\d(w_1), \\     %
\d(w_1)\d(w_2) %
=\d(w_2)\d(w_1)     %
\endgather           %
$$
where $\al, \bt, \gm, \dl$ are functions of $w=w_1-w_2, \h,
 \la, \eta, \tau$.

An  {\it 
\A/-module structure} on a vector space $M$ over $\Fun/(\C)$
is 
a direct sum decomposition 
$$
M\,=\, \oplus_{\mu \in \C/{1\over \eta}\Z} M[\mu]
\tag 5
$$
and endomorphisms  $ \a(w),$ $ \b(w),$
$\c(w),$  $ \d(w) \in \End/_\C(M)$ which meromorphically depend on
the parameter $ w \in \C$.
The direct sum decomposition allows us to define
endomorphisms $f(\la,\h) \in \End/(M)$ by the rule 
$f(\la,\h)m=f(\la,\mu)m$, if
$m\in M[\mu]$. We assume that the endomorphisms $f(\la,\h)$,
$ \a(w),$ $ \b(w),$
$\c(w),$  $ \d(w)$ satisfy the relations of the operator algebra.

An \A/-module $M$ is of {\it finite type}, if 
each space $M[\mu]$ has finite dimension
over $\Fun/(\C)$.

An \E/-module structure on
$V$ determines an \A/-module structure on $\Fun/(V)$.
The converse is also true.

For a finite type \A/-module $M$, we construct an \E/-module $V$.  Let
$e_j[\mu], j\in J_\mu$, be a basis of $M[\mu]$ over $\Fun/(\C)$.
Let $V[\mu]$ be a complex vector space with a basis denoted by
$e'_j[\mu], j\in J_\mu$. Set $V\,=\, \oplus_\mu V[\mu]$. 
Define an \E/-action on $V$. Set $f(h)\,e_j'[\mu]\,=\, f(\mu)\,e_j'[\mu]$.
To define an action of the other operators write
$$
\gather
\a(w)\,e_j[\mu]\,=\, \sum_k A_j^k(\la,w)\,e_k[\mu],
\qquad \b(w)\,e_j[\mu]\,=\, \sum_k B_j^k(\la,w)\,e_k[\mu-2], \\
\c(w)\,e_j[\mu]\,=\, \sum_k C_j^k(\la,w)\,e_k[\mu+2], 
\qquad \d(w)\,e_j[\mu]\,=\, \sum_k D_j^k(\la,w)\,e_k[\mu],
\endgather
$$
for suitable  functions
$A^k_j(\la,w), B^k_j(\la,w), C^k_j(\la,w), D^k_j(\la,w) $ $\in \Fun/(\C)$
and then set
$$
\gather
a(\la,w)\,e'_j[\mu]\,=\, \sum_k A_j^k(\la,w)\,e'_k[\mu],
\qquad b(\la,w)\,e'_j[\mu]\,=\, \sum_k B_j^k(\la,w)\,e'_k[\mu-2], \\
c(\la,w)\,e'_j[\mu]\,=\, \sum_k C_j^k(\la,w)\,e'_k[\mu+2], 
\qquad d(\la,w)\,e'_j[\mu]\,=\, \sum_k D_j^k(\la,w)\,e'_k[\mu],
\endgather
$$
These formulae define an \E/-module structure on $V$.

Let $V$ and $W$ be complex vector spaces.
Any meromorphic $1$-periodic
function $\phi (\la)$ with values in $\Hom/_\C(V,W)$
induces a homomorphism $\Fun/(V) \to \Fun/(W),$
$ F(\la) $ $\mapsto $ $
\phi (\la)F(\la).$

A {\it morphism} of an \E/-module $V$ to an \E/-module $W$ is
a $1$-periodic meromorphic function $\phi (\la)$ with values in $\Hom/_\C(V,W)$
such that the induced homomorphism 
$\Fun/(V) \to \Fun/(W)$ commutes with the action of the
operators 
$ f(\la,\h), \, \a(w),$ $ \b(w), \,$ $ \c(w),$ $  \, \d(w)$. 
A morphism is  an {\it isomorphism},  if the homomorphism
$\phi(\la) $ is nondegenerate for generic $\la$.

If $\phi_1(\la) \in \Hom/(V_1, W_1)$ and $\phi_2(\la) \in \Hom/(V_2, W_2)$
are morphisms, then
$$
(\phi_1 \otimes \phi_2) (\la) \,=\,\phi_1^{(1)}(\la - 2\eta h^{(2)})
   \phi_2^{(2)} (\la) 
$$
is a morphism from $V_1 \otimes V_2$ to $W_1 \otimes W_2$.

An \E/-module $W$ is {\it irreducible}, if for all non-trivial morphisms
$\phi(\la) : V \to W$ the map $\phi(\la)$ is surjective for generic $\la$.
A module is {\it reducible}, if it is not irreducible.

A {\it singular vector} in an \E/-module $V$ is a non-zero element
$v\in \Fun/(V)$ such that $\c(w)v\, =\, 0$  for all $w$.
An element $v\in \Fun/(V)$ is of {\it $h$-weight} $\mu$, if
$f(\la,\h)v=f(\la,\mu)v$ for all $f(\la,\h)$.
An element $v\in \Fun/(V)$ is of {\it weight} $(\mu, A(\la,w), 
D(\la,w))$, if it is of $h-$weight $\mu$ and
$\a(w)v = A(\la,w)v$, $\d(w)v = D(\la,w)v$  for all $w$.

Let $A(\la,w)$ and $D(\la,w)$ be two functions, $\La \in \Ce/$ .
An \E/-module is a {\it highest weight module} with highest weight
$(\La, A(\la,w), D(\la,w))$ and highest weight vector $v\in 
\Fun/(V)$, if $v$ is a singular vector of weight
$(\La,$ $ A(\la,w),$ $ D(\la,w))$ and 
 if $\Fun/(V)$ is generated over $\Fun/(\C)$ by elements of the form
$\b(w_1)\cdot ...\cdot \b(w_n)\,v$, where $w_1,...,w_n$ is an arbitrary
finite set.

In this paper we will consider 
only the highest weight modules with highest weights
$(\La,$ $ A(\la,w),$ $ D(\la,w))$  such that the functions
$A(\la,w)$ and $D(\la,w)$ are not identically equal to zero.

\proclaim{Theorem 1}

A highest weight \E/-module is  reducible, if and only if
the module has a singular vector of $h$-weight $\La-2k$ for a positive
integer $k$.

\endproclaim

\proclaim{Theorem 2}

Let $V$ and $W$ be irreducible highest weight \E/-modules with highest weight
vectors $v \in \Fun/(V)$ and $w \in \Fun/(W)$ and the same highest weight.
Then any isomorphism $V \to W$ sends the highest weight vector
of the first module to the highest weight vector of the second module multiplied
by a scalar meromorphic  doubly periodic function $g(\la)$ with periods $1$ and
$2\eta$. Moreover, for any such a function $g(\la)$ there exists a unique
isomorphism $V \to W$ sending $v$ to $g(\la)\,w$.

\endproclaim

\head
4. Evaluation modules
\endhead

For  complex numbers $\La$ and $z$ we describe an infinite
dimensional module $V_{\La}(z)$
of the \eqg/ which we call the {\it evaluation Verma module}.
If $\La=n+(m+\ell\tau)/2\eta,$
 where $ n, m, \ell$ are integers, $n \geq 0$,
then the module has a submodule such that the quotient module
$L_{\La}(z)$ has finite dimension $n+1$.

Let $V_{\La}(z)$ be an infinite dimensional complex vector space with a basis
$e_k,\, k\in \Z_{\ge 0}$.  Define an action of $f(h)$ by
$$
f(h)\,e_k\,=\,f(\La-2k)e_k.      
$$
and an action of the other generators  by
$$
\gather %
a(\la,w)\,e_k\,=\,{\theta(z-w+(\La+1-2k)\eta)\over
\theta(z-w+(\La+1)\eta)} {
\theta(\la+2k\eta) \over 
\theta(\la) }\, e_k,              \\
b(\la,w)\,e_k\,=\,-\,
{\theta(-\la+z-w+(\La-1-2k)\eta)\over
\theta(z-w+(\La+1)\eta)}{             
\theta(2\eta) \over \theta(\la) }\, e_{k+1}, \\             
c(\la,w)\,e_k\,=\,-\,{\theta(-\la-z+w+(\La+1-2k)\eta)\over
\theta(z-w+(\La+1)\eta)}    {         
\theta(2(\La+1-k)\eta)\over
\theta(\la) } {\theta(2k\eta)\over \theta(2\eta)} \,e_{k-1} \\
d(\la,w)\,e_k\,=\,{\theta(z-w+(-\La+1+2k)\eta)\over
\theta(z-w+(\La+1)\eta) }{              
\theta(\la-2(\La-k)\eta) \over 
\theta(\la) }\, e_k.\\              
\endgather
$$

\proclaim{Theorem 3}

These formulae define an \E/-module structure on
$V_\La(z)$. If
 $\La=n + (m + \ell\tau)/2\eta$, where $ n, m, \ell$ are integers,
 $n \geq 0$,
then the subspace spanned by  $e_k$, $k > n$, is a submodule.
The quotient space
$L_\La(z)$ is a module of  dimension $n+1$. 

\endproclaim

The module  is a highest weight module with highest weight vector
$e_0$ and
highest weight
$(\La, A(\la,w), D(\la,w))$, where 
$A(\la,w)\,=\,1$ and the function 
$D(\la,w))$ is determined by $\La$ and $z$, namely,
$$
D(\la,w)\,=\,{\theta(z-w+(-\La+1)\eta)\over
\theta(z-w+(\La+1)\eta) }{              
\theta(\la-2\La\eta) \over 
\theta(\la) }\,.
$$
Here $e_0$ is considered as  a constant function in $\Fun/(V_\La(z))$.

\proclaim{Remark}

Let $z(\la)$ be a doubly periodic meromorphic function with periods
$1$ and $2\eta$. Let $\La$ be a complex number. Then
the above formulae define an \E/-module structure on a complex space
$V_{\La}(z(\la))$ with a basis $e_k$, if we substitute $z=z(\la)$.
 If
$\La=n + (m + \ell\tau)/2\eta$, where $ n, m, \ell$ are integers,
 $n \geq 0$,
then the subspace spanned by  $e_k$, $k > n$, is a submodule.
The quotient space
$L_\La(z(\la))$ is a module of  dimension $n+1$. 
\endproclaim

\proclaim{Theorem 4}

For generic $\eta$ the evaluation module $V_\La(z)$
is 
reducible if and only if  $\La=n + (m + \ell\tau)/2\eta$, where $ n, m, \ell$ 
are integers,  $n \geq 0$.

For generic $\eta$ the finite 
dimensional module  $L_{n + (m + \ell\tau)/2\eta}(z)$
is irreducible.

\endproclaim

Denote by $W_\La(z)$ the finite dimensional
module $L_\La(z)$, if  $\La=n + (m + \ell\tau)/2\eta$
where $m, \ell$ are integers,  $n \geq 0$, and the infinite dimensional
module $V_\La(z)$ if  $\La$ does not have this form.

\proclaim{Theorem 5}

Let $\eta $ be generic, $\La_1,...,\La_n$ arbitrary. Then for generic
$z_1,...,z_n$ the module
$$
W_{\La_1,...,\La_n}(z_1,...,z_n)\,= \,
W_{\La_{1}}(z_{1}) \otimes ... \otimes
W_{\La_{n}}(z_{n})
$$  
is an irreducible highest weight \E/-module
with highest weight \newline $(\La_1+...+\La_n,$ $ 1,$ 
$ D(z_1,...,z_n,\La_1,...,\La_n))$ ,
$$
D(z_1,...,z_n,\La_1,...,\La_n))\,=\, 
{\theta(\la - 2(\La_1+...+\La_n)\eta) \over \theta(\la)        }
\prod^n_{k=1} { \theta(z_k-w+(-\La_k+1)\eta) \over  
\theta(z_k-w+(\La_k+1)\eta) }    \, ,   
$$
and highest weight vector $e_0(1)\otimes ...\otimes e_0(n)$,
where $e_0(j)$ is the highest weight vector of the $j$-th
factor.

\endproclaim

\proclaim{Corollary of Theorems 2 and 5} 

 Two irreducible tensor products of  evaluation
modules are isomorphic, if the highest
weights of the products are equal (while their factors could 
have different $z_k$ and $\La_k$). 

\endproclaim

\head
5. One dimensional modules
\endhead

\proclaim {Theorem 6}

If an \E/-module is one dimensional, then $b(\la,w)$ and
$c(\la,w)$ act by zero for all $\la, w$ and the module 
is a highest weight module
with highest weight $(\La, A, D)$. The highest weight is of one
of the following three types.

If $A$ and $D$ are nonzero, then 
$2\eta\La = \ell\tau$ where $ \ell$ is an integer, 
$$
A(\la,w)\,=\, f(\la,w)\,g(\la), \qquad
D(\la,w)\,=\,e^{2\pi i\ell(\la-w)}\, f(\la,w)\,j(\la)/g(\la+2\eta),
$$
where $f(\la,w),\,g(\la),\,j(\la)$ are  arbitrary
meromorphic
functions such that
$f(\la,w)$ and $j(\la)$ are $2\eta$-periodic in $\la$ and
$A(\la,w), D(\la,w)$ are $1$-periodic in $\la$.

If $A \equiv 0$, then $\La$ is arbitrary and 
$D(\la,w)\,=\, f(\la,w)\,g(\la)$
where $D(\la,w),\,f(\la,w),\,g(\la)$ have the same periodicity
properties as before.

If $D \equiv 0$, then $\La$ is arbitrary and 
$A(\la,w)\,=\, f(\la,w)\,g(\la)$
where $A(\la,w),\,f(\la,w),\,g(\la)$ have the same periodicity
properties as before.

\endproclaim

Denote the module corresponding to nonzero $A$ and $D$
 by $U_{\ell\tau,f,g,j}$ and by $U_{\ell\tau,j}$,
if $f = g = 1$.

An important special case is formed by the modules with
$f = g = 1, \, j =$ const.

We have
$$
U_{\ell_1\tau,f_1,g_1,j_1} \otimes U_{\ell_2\tau,f_2,g_2,j_2}\,=\,
U_{(\ell_1+\ell_2)\tau,f,g,j} 
$$ 
where
$$
f(\la,w)=f_1(\la-\ell_2\tau,w)f_2(\la,w), \,\qquad 
g(\la)=g_1(\la-\ell_2\tau)g_2(\la), 
$$
$$
j(\la)=e^{-2\pi i \ell_1\ell_2\tau}j_1(\la-\ell_2\tau)j_2(\la).
$$

The set of one dimensional modules forms a group 
with respect to the tensor product,
$U_{0,1}$ is the unit element.

\head
6. Isomorphism of evaluation modules
\endhead

\proclaim {Theorem 7}

For any complex numbers $\La, z$ and an integer $m$, the evaluation modules
$ V_{\La+m/\eta}(z)$ and $V_\La(z) $
are isomorphic with an isomorphism
map $\phi(\la) : 
V_{\La+m/\eta}(z) \to V_\La(z),$ $ e_k \mapsto e_k$.

For any $\La, z$ and an integer $\ell$, the evaluation modules
$V_{\La+\ell\tau/\eta}(z)$  and \newline
$U_{\ell\tau,j_L}\otimes V_{\La}(z)
\otimes U_{\ell\tau,j_R}$ 
are isomorphic with an isomorphism
map 
$$
\phi(\la) : V_{\La+\ell\tau/\eta}(z) \to
U_{\ell\tau,j_L}\otimes V_{\La}(z)
\otimes U_{\ell\tau,j_R},\qquad \, e_k \mapsto v_L \otimes e_k \otimes v_R,
$$
where
$$
j_L\,=\,e^{2\pi i (z+(\La+1)\eta)\ell},
\qquad
j_R\,=\,e^{2\pi i (z+3(\La-1)\eta)\ell+4\pi i \ell^2\tau}.
$$
and $v_L, v_R$ are  generating vectors of the one
dimensional modules.

\endproclaim

\proclaim{Remark} \endproclaim

Theorems 3 and 7 show that for every non-negative integer $n$
there are at most four potentially non-isomorphic $n+1$-dimensional
evaluation modules modulo tensoring with one dimensional 
representations.

\head
7. Highest weights of finite type modules 
\endhead

In this section we discuss highest weights of highest weight 
\E/-modules of finite type. 
Recall that we consider only the modules with highest weights 
$(\La, \,A(\la,w),\, D(\la,w))$ such that both functions $A$ and $D$ are 
not identically equal to zero. 

\proclaim {Lemma 8}

If $(\La, A(\la,w), B(\la,w))$ is a highest weight, then
$$
A= f(\la,w)\,g(\la), \qquad 
D(\la,w)\,=\,{\theta (\la-2\eta\La)
\over \theta (\la)}\,
{h(\la,w) \over g(\la+2\eta)},
$$
where $f(\la,w),\,g(\la),\,h(\la,w)$ are  meromorphic
functions such that
$f(\la,w)$ and $h(\la,w)$ are $2\eta$-periodic in $\la$ and
$A(\la,w)$ and $ D(\la,w)$ are $1$-periodic in $\la$.

\endproclaim

Considering highest weights we could restrict ourself to the case
$A(\la,w)=1$. In fact, let
 $V$ be a highest weight \E/-module of finite type
with highest weight
$(\La,\, A(\la,w),\, D(\la,w))$. There is  a one dimensional
\E/-module $U$ with highest weight 
$(0,\, 1/A(\la,w),\, C(\la,w))$ with suitable $C(\la,w)$. 
Then $V\otimes U$
is a highest weight module of finite type with
highest weight $(\La,\, 1,\, D(\la,w)\,C(\la,w))$.

\proclaim {Theorem 9}

Let $\eta$ be real and 
 irrational.  Let a highest weight
\E/-module be of finite type and have highest weight of the form
$(\La, 1, D(\la,w))$ with nonzero function $D(\la,w)$. Then 
$D(\la,w)$ has the form
$$
D(\la,w)\,=\, F\,
{\theta(\la - 2\La\eta) \over \theta(\la)        }
\prod^n_{k=1} { \theta(s_k-w) \over  
\theta(t_k-w) }    \, ,   
\tag 6
$$
where $n\geq 0$ and $F,\, s_k,\, t_k$ are suitable constants.

\endproclaim

Note, that after tensoring the highest weight module
with a suitable one dimensional module,
 we can set $F=1$ in (6).

Consider the submodule of a tensor product of evaluation modules
generated by the tensor product of highest weight vectors of factors.
Consider the quotient of the submodule by its 
maximal proper submodule. The resulting irreducible highest weight
module will be called the {\it irreducible module associated with
the tensor product}.

\proclaim {Corollary }

Let $\eta$ be real and irrational, 
then any irreducible highest weight module of finite type with highest weight
$(\La,\, 1,\, D(\la,w))$, where $D(\la,w))$ is given by (6)
with $F=1$,  is isomorphic to the irreducible highest weight module
associated with a tensor
product of evaluation representations.

\endproclaim

\head
8. Fusion of evaluation modules
\endhead

\proclaim{Theorem 10}

The tensor product of two evaluation 
modules $V_{\La_1}(z_1)\otimes V_{\La_2}(z_2)$
contains the evaluation submodule $V_{\La_1+\La_2}(z_2-\La_1\eta)$,
if $z_2-z_1=(\La_1+\La_2)\eta$. The imbedding  
$V_{\La_1+\La_2}(z_2-\La_1\eta) \hookrightarrow 
V_{\La_1}(z_1)\otimes V_{\La_2}(z_1+(\La_1+\La_2)\eta)$
is defined by the following formula, involving elliptic binomial coefficients,
$$
e_j \, \mapsto \, \sum^j_{\ell=0}
{ \theta(2j\eta)\,\theta(2(j-1)\eta)\cdots \theta(2(\ell+1)\eta)
  \over  
\theta(2\eta)\cdots \theta(2(j-\ell)\eta)} \, e_\ell \otimes e_{j-\ell}\,.
\tag 7
$$

\endproclaim

\proclaim{Theorem 11}

Let $\eta$ be generic.
Let $\La_1, \La_2$ have the form $\La_k=n_k+(m_k+\ell_k\tau)/2\eta$
where $k=1,2$ and $n_k, m_k, \ell_k$ are integers, $n_k\geq 0$.  Let 
$V\,=\,L_{\La_1}(z_1) \otimes L_{\La_2}(z_2)$ be the tensor product
of two finite dimensional irreducible evaluation modules.
The tensor product is reducible if and only if
$$
z_1-z_2 \,\, {\text{or}} \,\, z_2-z_1 \, = \ (\La_1 + \La_2 - 2j + 2)\eta +
m + \ell\tau \,
$$
for some integers $j,\, \ell, m,\, 0 < j \leq {\text{min}}\{n_1, n_2\}.$
In this case the tensor product has a unique proper submodule $W$. Moreover:

If $z_2-z_1\,=\,(\La_1 + \La_2 - 2j + 2)\eta +
m + \ell\tau $, then
$$
\gather
W\, \simeq \, L_{j-1}(z_1+(\La_1-j+1)\eta) \otimes
L_{\La_1+\La_2-j+1}(z_2-(\La_1-j+1)\eta) \otimes U_1\,, \\
V\,/\,W\, \simeq \, L_{\La_1-j}(z_1-j\eta) \otimes
L_{\La_2-j}(z_2+j\eta) \otimes U_2\,, \\
\endgather
$$
where $U_1, U_2$ are suitable one dimensional modules.

If $z_1-z_2\,=\,(\La_1 + \La_2 - 2j + 2)\eta +
m + \ell\tau $, then
$$
\gather
W\, \simeq \, L_{\La_1-j}(z_1+j\eta) \otimes
L_{\La_2-j}(z_2-j\eta) \otimes U_1\,, \\
V\,/\,W\, \simeq \, L_{j-1}(z_1-(\La_1-j+1)\eta) \otimes
L_{\La_1+\La_2-j+1}(z_2+(\La_1-j+1)\eta) \otimes U_2\,,\\
\endgather
$$
where $U_1, U_2$ are suitable one dimensional modules.

\endproclaim

If $z_1-z_2\,=\,(\La_1 + \La_2 - 2j + 2)\eta +
m + \ell\tau $, then the imbedding $W\, \hookrightarrow \,V$
sends the highest weight vector of $W$ to a singular vector
of $V$ of $h$-weight $\La_1+\La_2-2j$. The singular vector is unique
up to multiplication by a scalar function of $\la$. The singular vector
has the form $\sum_\ell (-1)^\ell\,
A_\ell(\la)\, e_\ell \,\otimes\,e_{j-\ell}$ where 
$$
\gather
A_\ell(\la)\,=\, {\theta(2j\eta)\,\theta(2(j-1)\eta)
\cdots \theta(2(\ell +1)\eta)
\over \theta (2\eta) \theta (4\eta) \cdots \theta (2(j-\ell )\eta)}\, 
\prod^{\ell -1}_{i=0} 
{\theta(-\la+2\eta(\La_1+\La_2-2j+\ell -i+1))
\over \theta(-\la+2\eta(-j+\ell -i))} \times \\
\hfill \biggl[ \, \prod^{\ell -1}_{i=0} \theta(2(\La_1-i)\eta)) 
\cdot
\prod^{j-\ell -1}_{i=0} \theta(2(\La_2-i)\eta))\,\biggr]^{-1}.\\ 
\endgather
$$

If $z_2-z_1=(\La_1+\La_2)\eta$, then 
$W \simeq  L_{\La_1+\La_2}(z_2-\La_1\eta)$,
and the imbedding  
$L_{\La_1+\La_2}(z_2-\La_1\eta) \hookrightarrow 
L_{\La_1}(z_1)\otimes L_{\La_2}(z_1+(\La_1+\La_2)\eta)$
is given by formula (7).

\head 9. The universal evaluation module
\endhead

The formulae defining evaluation representations
admit the following generalization. 
Let $\V$ be the
complex vector space of all
functions in $\la$, $h$, $z$, $\La$.
 Define an action of the operator algebra \A/ on $\V$.
Namely, define an action of operators $f(\la,\h)$
by the rule
$$
(f(\la,\h)\, v)(\la,h,z,\La)\,=\,f(\la,h)\, v(\la,h,z,\La),
$$
where $v(\la,h,z,\La) \in \V$,
and an action of the other operators by the rule
$$
\gather
(\a(w)\,v)(\la,h,z,\La) \,=\,{\theta(z-w+(h+1)\eta)\over
\theta(z-w+(\La+1)\eta)} {
\theta(\la-(h-\La)\eta) \over 
\theta(\la) }\, v(\la-2\eta,h,z,\La),              \\
(\b(w)\,v)(\la,h,z,\La)\,=\,-\,
{\theta(-\la+z-w+(h-1)\eta)\over
\theta(z-w+(\La+1)\eta)}{             
\theta(2\eta) \over \theta(\la) }\, v(\la+2\eta,h-2,z,\La), 
\\             
(\c(w)\,v)(\la,h,z,\La)\,=\,-\,{\theta(-\la-z+w+(h+1)\eta)\over
\theta(z-w+(\La+1)\eta)}    {         
\theta((h+\La+2)\eta)\over
\theta(\la) } {\theta((\La-h)\eta)\over \theta(2\eta)} \,
v(\la-2\eta,h+2,z,\La), \\
(\d(w)\,v)(\la,h,z,\La)\,=\,{\theta(z-w+(-h+1)\eta)\over
\theta(z-w+(\La+1)\eta) }{              
\theta(\la-(h+\La)\eta) \over 
\theta(\la) }\, v(\la+2\eta,h,z,\La).\\              
\endgather
$$

\proclaim {Theorem 12}

The operators $f(\la,\h)$, $\a(w)$, $\b(w)$, $\c(w)$, $\d(w)$ $\in
\End/(\V)$ satisfy the relations
of the elliptic operator algebra \A/.

\endproclaim

$\V$ is called the {\it universal generalized
evaluation \A/-module}.

Precisely speaking $\V$ is not a module in the sense defined in Section
3.

The universal evaluation module has many invariant subspaces constructed 
in the following way. Let $\C^4$ be the complex space with coordinates
$\la,$ $h,$ $z,$ $\La$, and $X \subset \C^4$ a subset invariant with 
respect to the following four transformations
$$
\gather
(\la, h, z, \La) \, \mapsto \, (\la \pm 2\eta, h, z, \La) \,, \\
(\la, h, z, \La) \, \mapsto \, (\la , h \pm 2 , z, \La) \,.
\endgather
$$
Then the subspace $\V(X) \subset \V$ of all functions
with support in $X$ is invariant
with respect to the \A/-action.

For example, the set $X$ of lines $z=$const, $\La=$const,
$h = \Xi - 2k$, where $\Xi$ is a fixed number and $k\in \Z$, leads
to the following {\it cyclic \E/-module} $V_{\La,\Xi}(z)$.

$V_{\La,\Xi}(z)$ is 
an infinite dimensional complex vector space with a basis
$e_k,\, k\in \Z$.  The action of $f(h)$ is defined by
$$
f(h)\,e_k\,=\,f(\Xi-2k)\,e_k.      
$$
and the action of the other generators is defined  by
$$                  
\gather
a(\la,w)\, e_k \,=\,{\theta(z-w+(\Xi-2k+1)\eta)\over
\theta(z-w+(\La+1)\eta)} {
\theta(\la-(\Xi-\La-2k)\eta) \over 
\theta(\la) }\, e_k,              \\
b(\la,w)\,e_k\,=\,-\,
{\theta(-\la+z-w+(\Xi-2k-1)\eta)\over
\theta(z-w+(\La+1)\eta)}{             
\theta(2\eta) \over \theta(\la) }\, e_{k+1}, 
\\             
c(\la,w)\,e_k\,=\,-\,{\theta(-\la-z+w+(\Xi-2k+1)\eta)\over
\theta(z-w+(\La+1)\eta)}    {         
\theta((\Xi+\La+2-2k)\eta)\over
\theta(\la) } {\theta((\La-\Xi+2k)\eta)\over \theta(2\eta)} \,
e_{k-1}, \\
d(\la,w)\,e_k\,=\,{\theta(z-w+(-\Xi+2k+1)\eta)\over
\theta(z-w+(\La+1)\eta) }{              
\theta(\la-(\Xi+\La-2k)\eta) \over 
\theta(\la) }\, e_k.\\              
\endgather
$$

If $\La=\Xi$, then the subspace of 
$V_{\La,\Xi}(z)$ spanned by the vectors $e_k,\, k \geq 0$, forms a submodule
which is the evaluation module $V_\La(z)$.

The module $V_{\La}(z(\la))$  described in Section 4
also could be obtained by this construction.

Theorem 12 is an easy corollary of Theorem 3.

\subhead Remark
\endsubhead
One can construct universal evaluation modules for the Yangian
$Y(sl_2)$ and the  quantum loop group
$ U_q(\widetilde {  sl   }_2)$ by similar formulae. It might be that
these modules are new.

\head 10. Determinant \endhead

The element 
$$
\gather
\Det (\la,w)\,=\, { \theta(\la)\over
\theta(\la-2\eta h)}\,
( d(\la,w+2\eta)a(\la+2\eta,w)
- b(\la,w+2\eta)c(\la+2\eta,w)) \,=\\
=\,{ \theta(\la)\over
\theta(\la-2\eta h)}\,( a(\la,w+2\eta)d(\la-2\eta,w)
- c(\la,w+2\eta)b(\la-2\eta,w))\,\\
\endgather
$$
is called the {\it determinant element} of the \eqg/ \E/,
the corresponding element 
$$
\gather
\Dett (w)\,=\, { \theta(\la)\over
\theta(\la-2\eta \h)}\,
( \d(w+2\eta)\a(w)
- \b(w+2\eta)\c(w)) \,=\\
=\,{ \theta(\la)\over
\theta(\la-2\eta \h)}\,( \a(w+2\eta)\d(w)
- \c(w+2\eta)\b(w))\,\\
\endgather
$$
is called the {\it determinant element} of the elliptic quantum algebra \A/.

\proclaim {Theorem 13}

The determinant element $\Dett (w)$ is a central element in the elliptic
quantum algebra \A/.

\endproclaim

\proclaim {Theorem 14}

The determinant element $\Det (\la,w)$ is a group-like element in the elliptic
quantum group \E/. Namely, if $V$ and $W$ are \E/-modules, then 
$\Det (\la,w)$ acts in the module $V\otimes W$ as
$$
\Det (\la-2\eta h^{(2)},w) \otimes \Det (\la,w).
$$

\endproclaim

Note that other formulae for the determinant element could be deduced
>from  relations (4).

\head 11. Dual modules \endhead

Define a homomorphism  \E/ $ \to\, \C$, the {\it counit},
by the rule $f(h) \, \mapsto \, f(0)$,
$$
\gather  
a(\la,w)\,\mapsto \,1, \qquad b(\la,w)\,\mapsto \,0,\\
c(\la,w)\,\mapsto \,0, \qquad d(\la,w)\,\mapsto \,1.
\endgather
$$
The counit defines an \E/-module structure on $\C$.

Let
\E/ $\to \End/(V)$ be an \E/-module structure on
a complex vector space $V$. Assume that the determinant
element $\Det (\la,w) \in \End/(V)$ is nondegenerate for
generic $\la$ and $ w$ and denote by $\Det ^{-1}(\la,w)$
its inverse. 
Let
$$
V^*\,=\,\oplus_\mu \,V[\mu]^*
$$
be the restricted dual space to
  $V$. We introduce an \E/-module structure
on the restricted dual. 
The restricted dual space
 with this \E/-module structure is
called the module {\it dual} to $V$.

Introduce linear maps $Sf(h),\,Sa(\la,w),\,$ $Sb(\la,w),\,$ $Sc(\la,w),\,
Sd(\la,w)\, \in\,\End/(V)$
by the rule  
$$
\gather  
Sf(h) \, =\, f(-h)\, \tag 8 \\
Sa(\la,w)\,= \, {\theta(\la + 2\eta h - 2\eta) \over
\theta(\la -2\eta)}\,
 \Det ^{-1}(\la + 2\eta h - 2\eta, w)\,
d(\la + 2\eta h - 2\eta, w + 2\eta)\,, \\
Sb(\la,w)\,= \, -\, {\theta(\la + 2\eta h + 2\eta) \over
\theta(\la + 2\eta )}\,
\Det ^{-1}(\la + 2\eta h + 2\eta, w) \,
b(\la + 2\eta h - 2\eta, w + 2\eta)\,, \\
Sc(\la,w)\,= \,-\, {\theta(\la + 2\eta h - 2\eta) \over
\theta(\la -2\eta)}\,
\Det ^{-1} (\la + 2\eta h - 2\eta, w)\,
c(\la + 2\eta h + 2\eta, w + 2\eta)\,, \\
Sd(\la,w)\,= \,  {\theta(\la + 2\eta h + 2\eta) \over
\theta(\la + 2\eta )}\,
\Det ^{-1}(\la + 2\eta h + 2\eta, w)\,
a(\la + 2\eta h + 2\eta, w + 2\eta)\,. \\
\endgather
$$
Let
$Sf(h)^*,\,Sa(\la,w)^*,\,$ $Sb(\la,w)^*,\,$ $Sc(\la,w)^*,\,
Sd(\la,w)^*\, \in\,\End/(V^*)$ be their dual maps, respectively.

\proclaim {Theorem 15}

Let $V$ be an \E/-module. Then
the linear maps $Sf(h)^*,\,Sa(\la,w)^*,\,$ $Sb(\la,w)^*,\,$ $Sc(\la,w)^*,\,
Sd(\la,w)^*\, \in\,\End/(V^*)$ define an \E/-module structure on $V^*$.

Consider  $V\otimes V^*$ and  $ V^* \otimes V$ as \E/-modules. Then 
the natural maps
 $\C \, \to  V\otimes V^*$ and  $ V^* \otimes V \to \C$ are homomorphisms
of \E/-modules.

\endproclaim

Note that $\Det (\la,w)$ acts in $V^*$ as $\Det ^{-1}(\la + 2\eta h, w)^*$.

Formulae (8) play the role of the {\it antipode} for the \eqg/ \E/.

Repeating the above construction we introduce a new module
structure on the initial vector space $V= (V^*)^*$. In this case the generators
$f(h),\,a(\la,w),\,$ $b(\la,w),\,$ $c(\la,w),\,
d(\la,w)$ act as $f(h)$,
$$
\gather
{\theta(\la - 2\eta h - 2\eta) \over
\theta(\la - 2\eta h)}\,
{\theta(\la) \over
\theta(\la - 2\eta)}\,
 {\Det (\la, w)\over
\Det (\la, w + 2\eta)}\,
a(\la, w + 4\eta)\,, \\
{\theta(\la - 2\eta h - 2\eta) \over
\theta(\la - 2\eta h )}\,
{\theta(\la ) \over
\theta(\la +2\eta)}\,
 {\Det (\la, w)\over
\Det (\la, w + 2\eta)}\,
b(\la, w + 4\eta)\,, \\
{\theta(\la - 2\eta h + 2\eta) \over
\theta(\la - 2\eta h )}\,
{\theta(\la ) \over
\theta(\la - 2\eta)}\,
 {\Det (\la, w)\over
\Det (\la, w + 2\eta)}\,
c(\la, w + 4\eta)\,, \\
{\theta(\la - 2\eta h +2\eta ) \over
\theta(\la - 2\eta h )}\,
{\theta(\la ) \over
\theta(\la + 2\eta)}\,
 {\Det (\la, w)\over
\Det (\la, w + 2\eta)}\,
d(\la, w + 4\eta)\,, \\
\endgather
$$
respectively, and the determinant element $\Det (\la,w)$ acts
as $\Det (\la,w)$.

\subhead Remark
\endsubhead
These formulae inspire the following two constructions of
automorphisms of the elliptic quantum group. Namely,
for any $1$-periodic nonzero meromorphic function
$g(\la)$ one can define two automorphisms $I, J\,:\, $\E/ $\to$ \E/ 
by
$$
\gather
I\,:\,f(h)\,\mapsto\,f(h) , \\
I\,:\,a(\la,w)\,\mapsto\, g(\la-2\eta)\,a(\la,w) , \qquad
I\,:\,b(\la,w)\,\mapsto\, g(\la)^{-1}\,b(\la,w) ,
  \\
I\,:\,c(\la,w)\,\mapsto\, g(\la-2\eta)\,c(\la,w) , \qquad
I\,:\,d(\la,w)\,\mapsto\, g(\la)^{-1}\,d(\la,w) ,\\
\endgather
$$
and
$$
\gather
J\,:\,f(h)\,\mapsto\,f(h) , \\
J\,:\,a(\la,w)\,\mapsto\, g(\la - 2\eta h -2\eta)^{-1}\,a(\la,w) , \qquad
J\,:\,b(\la,w)\,\mapsto\, g(\la - 2\eta h -2\eta)^{-1}\,b(\la,w) ,   \\
J\,:\,c(\la,w)\,\mapsto\, g(\la-2\eta h)\,c(\la,w) , \qquad
J\,:\,d(\la,w)\,\mapsto\, g(\la-2\eta h)\,d(\la,w) . \\
\endgather
$$
These automorphisms preserve the determinant element,
$I,J\,:\, \Det (\la,w) \mapsto \Det (\la,w)$.

\head 12. Weyl group \endhead

The map 
$$
\gather
 f(h) \, \mapsto \, f(-h), \\
a(\la,w)\,\mapsto \,d(-\la,w), \qquad b(\la,w)\,\mapsto \,c(-\la,w),\\
c(\la,w)\,\mapsto \,b(-\la,w), \qquad d(\la,w)\,\mapsto \,a(-\la,w)
\endgather
$$
defines an automorphism of order two of the \eqg/ \E/.

This fact follows from  identities $\al (-\la) = \dl (\la)$
and $\beta (-\la) = \gamma (\la)$, where  $\al,$ 
$\bt,$ $\gm,$ $\dl$ are defined in Section 2.

\head 13. Solutions to the modified Yang-Baxter equation \endhead

In this section we construct solutions to 
the modified Yang-Baxter equation.

Let $V_k$ be an irreducible highest weight module with
highest weight vector $v_k \in \Fun/(V_k)$ and highest weight
$(\La_k, 1, D_k(\la,w)), \, k=1, 2, 3.$
Assume that  all tensor products
$V_i \otimes V_j$ and $V_i \otimes V_j \otimes V_k$ are irreducible
highest weight modules
where $i, j, k$ are pair-wise distinct.
Fix highest weight vectors
$v_i(\la- 2\eta \La_j) \otimes v_j(\la)$ in $\Fun/(V_i \otimes V_j)$ 
and $v_i(\la- 2\eta (\La_j+\La_k)) \otimes v_j(\la- 2\eta \La_k)
\otimes v_k(\la)$ in $\Fun/(V_i \otimes V_j \otimes V_k)$. 

Assume that the highest weights
of the highest weight vectors 
$v_i(\la- 2\eta \La_j) \otimes v_j(\la)$ in $\Fun/(V_i \otimes V_j)$ 
and 
$v_j(\la- 2\eta \La_i) \otimes v_i(\la)$ in $\Fun/(V_j \otimes V_i)$ 
are the same and, moreover, the highest weights
of the highest weight vectors
$v_i(\la- 2\eta (\La_j+\La_k)) \otimes v_j(\la- 2\eta \La_k)
\otimes v_k(\la)$  do not depend on the order of the numbers
$i, j, k$.

Note that irreducible tensor products of evaluation modules
have all these properties, cf. Theorem 5. These properties also
hold, if $\La_j$ are integers, see Lemma 8.

Let $R^\vee_{V_jV_i}(\la)\,:\, V_j\otimes V_i \, \to 
V_i\otimes V_j$ be the unique isomorphism of the modules
$V_j\otimes V_i$ and
$V_i\otimes V_j$ sending the distinguished highest weight
vector to the distinguished highest weight vector.
Set $R_{V_iV_j}(\la) = R^\vee_{V_jV_i}(\la)\, P$
where $P\,:\, V_i\otimes V_j \, \to 
 V_j\otimes V_i$ is the permutation of factors.

\proclaim {Theorem  16  }

The linear operators
$R_{V_iV_j}(\la) \, \in
\, {\text{Hom}}(V_i,V_j)$ satisfy the modified Yang-Baxter relation,
$$
\gather
R_{V_1V_2}(\la-2\eta h^{(3)})\, R_{V_1V_3}(\la)\,
R_{V_2V_3}(\la-2\eta h^{(1)})\, =           \\
R_{V_2V_3}(\la)\,R_{V_1V_3}(\la-2\eta h^{(2)})\,
R_{V_1V_2}(\la)\, ,\\
\endgather
$$
and the relation
$$
R^{(12)}_{V_1V_2}(\la)\,R^{(21)}_{V_2V_1}(\la)\, =\, {\text{Id}}.
$$

\endproclaim

Let $\eta$ be generic. Let $W_{\Lambda'}(z')$
and $W_{\Lambda''}(z'')$ be two irreducible evaluation modules
such that
$W_{\Lambda'}(z')\otimes W_{\Lambda''}(z'')$
and
$W_{\Lambda''}(z'')\otimes W_{\Lambda'}(z')$
are irreducible. Then the unique isomorphism
$$
R^\vee(\lambda)\, :\, W_{\Lambda'}(z')\otimes W_{\Lambda''}(z'')\,
\to\, W_{\Lambda''}(z'')\otimes W_{\Lambda'}(z')
$$
sending $e_0'\otimes e_0''$
to $e_0''\otimes e_0'$ can be constructed rather explicitly.
Namely, introduce elements $\b^{[k]}(w)$
and $\b_{[k]}(w)$ 
of the operator algebra
by formulae
$$
\gather
\b^{[k]}(w)\,=\,\b(w-2(k-1)2\eta)\cdots \b(w-2\eta)\,\b(w), \\
\b_{[k]}(w)\,=\,\b(w+2(k-1)2\eta)\cdots \b(w+2\eta)\,\b(w) \\
\endgather
$$
For any $k\geq 0,\,p>0$ set 
$\b_{k,p}\,=\,\b^{[k]}(z'+(\Lambda'+1)\eta)\,
\b_{[p]}(z''+(-\Lambda''+1)\eta)$ and set
$$
\gather
\b_{k,0}\,=\,\b^{[k-1]}(z'+(\Lambda'-1)\eta)\,\cdot
(\b(z'+(\Lambda'-1)\eta) \,-\, \\
\hfill
{\theta (-\lambda +z''-z'+(\Lambda''-\Lambda')\eta)
\over 
\theta (z''-z'+(\Lambda''-\Lambda'+2)\eta)}\,
{\theta (2\Lambda''\eta)
\over 
\theta (-\lambda +2(\Lambda''-1)\eta)}
\b(z''+(-\Lambda''+1)\eta)).\\
\endgather
$$
Then for all $k,p$
we have
$\b_{k,p}\,(e_0'\otimes e_0'')\,=\,
f_{k,p}(\lambda, z'-z'')\,
e_k'\otimes e_p''$
where $f$ is a scalar function. The scalar function is equal
to an alternating product of theta functions which can be written 
explicitly using  formulae of Section 4. Then
$$
R^\vee (\lambda)\,:\,
e_k'\otimes e_p''\,\mapsto \,
{1 \over f_{k,p} (\lambda, z'-z'')}\,
\b_{k,p}\,
(e_0''\otimes e_0')\,.
$$

\head 14. Commuting elements \endhead

\proclaim {Theorem 17}

For an \A/-module $M$ the endomorphisms $\t(w)\,=\,\a(w) + \d(w)$
pair-wise commute for all $w$ on the zero $h$-weight subspace $M[0]$. 

\endproclaim
                                                                                
This fact is related to the integrability of
the interaction-round-a-face  models of statistical mechanics
connected with the \eqg/ \E/, see [JMO], [JKMO] and references threin.
 We will describe the Bethe ansatz for
these models in the next paper [FTV].

\head 15. The case $2N\eta \,=\, 1$
\endhead

Representation theory of the elliptic quantum group
becomes special, if the parameter $\eta$ has the form $2N\eta\,=\,
m + \ell \tau$ where $ N, \, m, \, \ell\,$ are integers.
We will discuss this subject in a separate paper. Here we
make remarks about the case $2N\eta \,=\, 1$ where $N$ is a natural
number.
First
we construct two families of $N$-dimensional \E/-modules
$T_{\La}(z)$ and $T_{\La,\Xi}(z)$ where $\La, \Xi \in \C/\Z$
and $z \in \C$, and then we indicate some central elements of the operator
algebra \A/.

Consider an evaluation module $V_\La(z)$.  If $2N\eta=1$, then
the subspace spanned by
$e_k,\, k \geq N$, is a submodule. The quotient space
$T_\La(z)$ is an $N$-dimensional  module.

Consider a cyclic module  $V_{\La,\Xi}(z)$ with its distinguished basis
$e_k,\, k \in \Z$. Introduce a new basis $v_k,\, k \in \Z$. Namely, if
 $N$ is odd, then set
$$
 v_k\,=\, (-1)^{k(k+1)/2}\, e_k,
$$
and if $N$ is even, then set
$$
 v_k\,=\, (-1)^{k(k+1)/4}\, e_k.
$$
With respect to the new basis, the formulae for the action of
all operators of the \eqg/ \E/ become $N$-periodic in $k$.
Identifying $v_k$ and $v_{k+N}$ we get an $N$-dimensional
\E/-module $T_{\La,\Xi}(z)$.

\vskip3ex

Let $2\eta = 1$. Then  the operator
algebra $A_{\tau, \eta = 1/2}(sl_2)$ is commutative.

There is the following generalization of this fact. Let $N$ be
 a natural number, $z$
a complex number. Introduce
elements $\a^{(N)}(z),\,$ $\b^{(N)}(z),\,$ $\c^{(N)}(z),\,$ $\d^{(N)}(z)\,$
of the operator algebra $A_{\tau, \eta}(sl_2)$ by the formulae
$$
\gather
\a^{(N)}(z)\,=\, \a(z)\,\a(z-2\eta)\,\a(z-4\eta)\cdots \a(z-2(N-1)\eta),\,
\\
\b^{(N)}(z)\,=\, \b(z)\,\b(z-2\eta)\,\b(z-4\eta)\cdots \b(z-2(N-1)\eta),\,\\
\c^{(N)}(z)\,=\, \c(z)\,\c(z-2\eta)\,\c(z-4\eta)\cdots \c(z-2(N-1)\eta),\,
\\
\d^{(N)}(z)\,=\, \d(z)\,\d(z-2\eta)\,\d(z-4\eta)\cdots \d(z-2(N-1)\eta)\, .\\
\endgather
$$
For any element of the operator algebra \A/ of the form
$f(\la,\h)$, where $f$ is a meromorphic function  $1$-periodic 
in $\la$ and $1/2\eta$-periodic in $\h$, introduce a new element
$f^{(N)}(\la,\h)$ by the formula $f^{(N)}(\la,\h) = f(N\la,N\h)$.

 Let $V(\eta)$ be an \A/-module meromorphically
depending on $\eta$ and such that for $\eta = 1/(2N)$ the module 
is well defined.

As an example of such a module we can consider the \A/-module associated to
a tensor product of evaluation or cyclic \E/-modules.

\proclaim {Theorem 18}

Let $2N\eta\,=\,1$.
Consider the action of the operator algebra
$A_{\tau, 1/2N}(sl_2)$ in the module
 $V( 1/2N)$. Then each   of the operators $f^{(N)}(\la,\h),\,$
$\a^{(N)}(z),\,$ $\b^{(N)}(z),\,$ $\c^{(N)}(z),\,$ $\d^{(N)}(z),\,z \in \C,$
commutes with all of the operators
of the operator algebra.

\endproclaim

The proof of Theorem 18 is based on the following two facts.

Let $R(\la, w, \eta, \tau ) \in \End/(\C^2\otimes \C^2)$ be the
$R$-matrix defined in the introduction, then
the image of the linear map
$R(\la, -2\eta, \eta, \tau )$ is three dimensional and is
generated by vectors $v_+\otimes v_+$,
$v_+\otimes v_- + v_-\otimes v_+$,
$v_-\otimes v_-$, where $v_+=(1,0)$ and $v_-=(0,1)$.

If $2N\eta\, = \,1$, then the vectors
$v_+\otimes v_+$ and $v_-\otimes v_-$ lie in the kernel
of the linear map
Res$_{w = - 2(N-1)\eta} R(\la, w, \eta, \tau )$.

Now the proof of Theorem 18 follows from analysis of two Yang-Baxter-type
equations. If $4\eta=1$, then these two equations are
$$
\gather
R^{(21)}(\la-2\eta h, -2\eta, \eta, \tau)\, L^{(2)}(\la, z-2\eta)\,
L^{(1)}(\la-2\eta h^{(2)}, z)\, =           \\
L^{(1)}(\la, z)\,L^{(2)}(\la-2\eta h^{(1)}, z-2\eta)\,
R^{(21)}(\la, -2\eta, \eta, \tau)\,             \\
\endgather
$$
and
$$
\gather
R^{(13)}(\la-2\eta h, z-w, \eta, \tau)\,
R^{(23)}(\la-2\eta (h+h^{(1)}), z-2\eta - w, \eta, \tau)\,\times \\
\hfill
 L^{(1)}(\la, z)\,
L^{(2)}(\la-2\eta h^{(1)}, z-2\eta)\,
L^{(3)}(\la-2\eta (h^{(1)} + h^{(2)}), w)\,
=           \\
=\,L^{(3)}(\la, w)\, L^{(1)}(\la -2\eta h^{(3)}, z)\,
L^{(2)}(\la-2\eta (h^{(1)}+h^{(3)}), z-2\eta)\,\times \\
\hfill
R^{(13)}(\la, z-w, \eta, \tau)\,
R^{(23)}(\la-2\eta h^{(1)}, z-2\eta - w, \eta, \tau)\,.
\\
\endgather
$$
If $2N\eta\,=\,1$, then the equations are similar.

\head 16. The extended \eqg/
\endhead

The {\it extended } \eqg/ $\Ex/ _{\tau,\eta}$ is
the algebra over $\C$
with generators of two types. The generators of the first type
are labelled by meromorphic functions $f(h)$
of one complex variable with period $1/\eta$.
The generators of the second type 
$a(\la,w),\, b(\la,w),
c(\la,w),\, d(\la,w)$  are labelled by complex numbers $\la, w \in
\C$ (without any periodicity requirements).
 The generators of the
algebra satisfy the relations (2)-(4), cf. Section 2.

\proclaim {Theorem 19}

The extended \eqg/s $\Ex/_{\tau ,\eta} $ and $\Ex/_{\tau + 1,\eta} $ 
are isomorphic with an isomorphism 
$\Ex/_{\tau ,\eta} \to \Ex/_{\tau + 1,\eta} $ 
given by
$f(h) \mapsto f(h)$,
$$
(a(\la,w), b(\la,w), c(\la,w), d(\la,w))\,\mapsto 
(a(\la,w), b(\la,w), c(\la,w), d(\la,w)).
$$

Let $\tau '=-1/\tau$. Then 
the extended \eqg/s  $\Ex/_{\tau , \eta} $ and $\Ex/_{\tau ', \eta \tau '} $ 
are isomorphic with an isomorphism 
$\Ex/_{\tau , \eta} \to \Ex/_{\tau ', \eta \tau '} $ given by
$f(h) \mapsto f(h)$,
$$
\gather
(a(\la,w),\,\, b(\la,w), \, \, 
 c(\la,w), \,\, d(\la,w)) \,\, \mapsto  \,\,
(a(\tau '\la, \tau 'w)\,e^{2\pi i \tau ' A(\la,h,w)}, \tag 9 \\
b( \tau ' \la, \tau ' w)\,e^{2\pi i \tau ' B(\la,h,w)}, \,\, 
c( \tau '\la, \tau ' w)\,e^{2\pi i \tau ' C(\la,h,w)}
,\,\, d( \tau ' \la, \tau ' w)e^{2\pi i \tau ' D(\la,h,w)})\,, \\ 
\endgather
$$
where $A,\,B,\,C,\,D$ are 
polynomials in $\la, h, w$,
$$
\gather
A\,=\,(h-\La)\,\eta \,(-w+h+3\eta), \\
B\,=\,w\,(\la + 2\eta-\eta h + \eta \La)\, - \,(\la + 2\eta) \,(h-1)\, \eta\,+\, 
\eta^2\,((h+1)^2-(\La+1)^2)/2,\\
C\,=\,w(-\la+(h+\La+2)\eta)+(h+1)(-\eta\la+\eta^2(h+\La+2))
+\eta^2(\La-h-2)(\La-h+2)/2, \\
D\,=\,\eta\,( h + \La )\,( w + \eta h - \la - \eta ).
\endgather
$$
Here $\La \in \C$ is a parameter of the isomorphism.

\endproclaim

Say that an \EE/-module $V$ 
is {\it modular} with parameter $\La \in \C$, if in this module
the operators $ a(\la,w),\, b(\la,w), \, 
 c(\la,w), \, d(\la,w)$ are $1$-periodic  with respect to $w$,
$$
\gather
a(\la,w+1)\,=\, a(\la,w), \qquad b(\la,w+1) \,=  
 b(\la,w), \\ 
c(\la,w+1)\,= \, c(\la,w),\qquad
 d(\la,w+1) \,= \, d(\la,w),\\ 
\endgather
$$
and, moreover,
$$
\gather
 a(\la,w+\tau)\,=\, a(\la,w)\, e^{2\pi i \eta ( h - \La )}\,, \\
b(\la,w+\tau)\,=\, b(\la,w)\, e^{2\pi i ( - \la + \eta ( h - \La - 2))}\,, \\
c(\la,w+\tau)\,=\, c(\la,w)\, e^{2\pi i ( \la - \eta ( h + \La + 2 ))}\,, \\
d(\la,w+\tau)\,=\, d(\la,w)\, e^{2\pi i \eta ( - h - \La )}\,, \\ 
 \endgather
$$
for some constant $\La \in \C$.

\subhead Example 
\endsubhead 
An evaluation module $V_{\La}(z)$ and a cyclic module $V_{\La,\Xi}(z)$
are modular with parameter $\La$.

If $V_k$ is modular with parameter $\La_k,\, k = 1, 2,$ then
$V_1 \otimes V_2$ is modular with parameter $\La_1 + \La_2$.

Let $\tau '=-1/\tau$. Let an $\Ex/_{\tau ', \eta \tau '} $-module
$V$ be modular with parameter $\La$. Consider the isomophism
$\Ex/_{\tau , \eta} \to \Ex/_{\tau ', \eta \tau '} $ given by (9)
with the same parameter $\La$.
This isomorphism induces an \EE/-module structure on $V$.

\proclaim {Theorem 20}

This \EE/-module structure on $V$ is modular with parameter $\La$.

\endproclaim

\head 17. Concluding remarks
\endhead

Consider the \eqg/ $E_{\tau,\eta}(\frak g)$
associated to a simple Lie algebra $\frak g$
of $A, B, C, D$
type. For this \eqg/
one defines the notions of the operator algebra, an irreducible module,
a highest weight module, a finite type module, a singular vector, 
the Weyl group
in the same way as for the elliptic quantum group associated to
$sl_2$. The basic interrelations between these notions (like
Theorems 1 and 2) remain true after obvious changes. We plan to describe
evaluation modules  of the \eqg/ $E_{\tau,\eta}(sl_N)$
in a future paper.

In [FTV] we will describe solutions to the quantum Knizhnik - Zamolodchikov
- Bernard equations with values in a tensor product of evaluation \E/-modules.

The proofs of all formulated Theorems  are straightforward and will be 
published elsewhere.

\head References
\endhead

\vskip3ex 

\roster

\item"[ABB]" 
J.Avan, O.Babelon, E.Billey, {\it The Gervais-Neveu-Felder
equation and the quantum Calogero-Moser systems, } preprint
hep-th/9505091.

\item"[BBB]" 
O.Babelon, D.Bernard, E.Billey, {\it A quasi-Hopf algebra interpretation
of quantum 3-j and 6-j symbols and difference equations
,} preprint q-alg/9511019

\item"[CP]" V.Chari, A.Pressley, {\it A guide to quantum groups,}
University Press, Cambridge, 1994.

\item"[Fe1]" G.Felder, {\it Conformal field theory and integrable
systems associated to elliptic curves,} preprint hep-th/9407154,
to appear in the Proceedings of the ICM, Zurich 1994.

\item"[Fe2]" G.Felder, {\it Elliptic quantum groups,} preprint  hep-th/9412207,
to appear in the Proceedings of the ICMP, Paris 1994.

\item"[FTV]" G.Felder, V.Tarasov, A.Varchenko,
 {\it Bethe ansatz for interaction-round-a-face
elliptic models and solutions to $q$KZB equations,} in preparation.

\item"[JKMO]"  M.Jimbo, A.Kuniba,
T.Miwa and M.Okado, 
 {\it The $A^{(1)}_n$ face models,}
Commun. Math. Phys. 119 (1988), 543-565.

\item"[JMO]" M.Jimbo,
T.Miwa and M.Okado, {\it  Solvable lattice models related to
vector representations of classical simple Lie algebras,}
Commun. Math. Phys. 116 (1988), 507-525.

\endroster

\vskip3ex 

Department of Mathematics, Phillips Hall, University of North Carolina at
Chapel Hill, Chapel Hill, NC 27599-3250, USA

{\it E-mail address: felder\@math.unc.edu, av\@math.unc.edu.}

\end

\head . Restricted algebra \endhead

The elliptic quantum group has many subalgebras. For
any complex number $\mu$ consider the {\it restricted elliptic quantum group}
 \Em/ which is the subalgebra of \E/
generated by all generators of the form $f(h)$, $a(\la,w)$,
$b(\la,w)$, $c(\la,w)$, $d(\la,w)$, where  $\la \in L(\mu) $ and
$$
L(\mu)\,=\,\{\la \,|\, \la = \mu + 2k\eta,\, k\in \Z.\}
$$
Note that each relation of the \eqg/ involves only $\la$'s having this form
for some  $\mu$.

An   \Em/-module structure on a complex vector space $V$ is
defined similarly to an \E/-module structure. 

Let $\Fun/(V,\mu)$ be the space of functions on $L(\mu)$ with values in
$V$. A {\it morphism} of an \Em/-module $V$ to an \E/-module $W$ is
a function $\phi (\la)$ on $L(\mu)$
with values in $\Hom/_\C(V,W)$
such that the induced homomorphism 
$\Fun/(V,\mu) \to \Fun/(W,\mu)$ commutes with the action of the
operators 
$ f(\la,\h), \, \a(w),$ $ \b(w), \,$ $ \c(w),$ $  \, \d(w)$ 
where the operators are defined by  formulae in Section 3. 
A morphism is  an {\it isomorphism},  if the homomorphism
$\phi(\la) $ is nondegenerate for ALL(?) $\la \in L(\mu)$.

An \E/-module $W$ is {\it irreducible}, if for all non-trivial morphisms
$\phi(\la) : V \to W$ the map $\phi(\la)$ is surjective for ALL(?) $\la$.
A module is {\it reducible}, if it is not irreducible.

AND SO ON ...

\end

A {\it singular vector} in an \E/-module $V$ is a non-zero element
$v\in \Fun/(V)$ such that $\c(w)v\, =\, 0$  for all $w$.
An element $v\in \Fun/(V)$ is of {\it $h$-weight} $\mu$, if
$f(\la,\h)v=f(\la,\mu)v$ for all $f(\la,\h)$.
An element $v\in \Fun/(V)$ is of weight $(\mu, A(\la,w), 
D(\la,w))$, if it is of $h-$weight $\mu$ and
$\a(w)v = A(\la,w)v$, $\d(w)v = D(\la,w)v$  for all $w$.

Let $A(\la,w)$ and $D(\la,w)$ be two functions, $\La \in \Ce/$ .
An \E/-module is a {\it highest weight module} with highest weight
$(\La, A(\la,w), D(\la,w))$ and highest weight vector $v\in 
\Fun/(V)$, if $v$ is a singular vector of weight
$(\La,$ $ A(\la,w),$ $ D(\la,w))$ and 
 if $\Fun/(V)$ is generated over $\Fun/(\C)$ by elements of the form
$\b(w_1)\cdot ...\cdot \b(w_n)\,v$, where $w_1,...,w_n$ is an arbitrary
finite set.

\end

\proclaim{Theorem 6 }

Let $\La_1,...,\La_n$ be arbitrary.  
Let $V_{\La_1,...,\La_n}(z_1,...,z_n)$ be an abstract
irreducible highest weight
\E/-module with highest weight
$(\La_1+...+\La_n,$ $ 1,$ $
D(z_1,...,z_n,\La_1,...,$ $\La_n))$, where
$D(z_1,...,z_n,$ $\La_1,...,\La_n)$  is given by the above formula.
Then for generic $z_1,$ $...,$ $z_n$ the module
$V_{\La_1,...,\La_n}(z_1,...,z_n)$ is isomorphic
to $W_{\La_1,...,\La_n}(z_1,...,z_n)$.

\endproclaim